\documentclass[12pt]{article}
\usepackage{amssymb}
\usepackage{latexsym,bm}
\usepackage{graphicx}
\usepackage{amsmath}
\usepackage{mathrsfs}

\newcommand{\bea}{\begin{eqnarray*}}
\newcommand{\eea}{\end{eqnarray*}}
\newcommand{\be}{\begin{equation}}
\newcommand{\ee}{\end{equation}}
\newcommand{\ben}{\begin{eqnarray*}}
\newcommand{\een}{\end{eqnarray*}}

\date{}
\begin{document}
\title{Pairs of a tree and a nontree graph with the same status sequence\footnote{E-mail addresses:
{\tt 235711gm@sina.com}(P.Qiao),
{\tt zhan@math.ecnu.edu.cn}(X.Zhan).}}
\author{Pu Qiao, Xingzhi Zhan\thanks{Corresponding author.}\\
{\small Department of Mathematics, East China Normal University, Shanghai 200241, China}
 } \maketitle
\begin{abstract}
 The status of a vertex $x$ in a graph is the sum of the distances between $x$ and all other vertices. Let $G$  be a connected graph.
 The status sequence of $G$ is the list of the statuses of all vertices arranged in nondecreasing order. $G$ is called status injective
 if all the statuses of  its vertices are distinct. Let $G$ be a member of a family of graphs $\mathscr{F}$ and 
 let the status sequence of $G$ be $s.$ $G$ is said to be status unique in $\mathscr{F}$ if $G$ is the unique graph in $\mathscr{F}$ whose 
 status sequence is $s.$ In 2011, J.L. Shang and C. Lin  posed the following two conjectures. Conjecture 1: A tree and a nontree graph cannot have the same status sequence. Conjecture 2:   Any status injective tree is status unique in all connected graphs. We settle these two conjectures negatively. For every integer $n\ge 10,$ we construct a tree $T_n$ and a unicyclic graph $U_n,$ both of order $n,$ with the following two properties:
 (1) $T_n$ and $U_n$ have the same status sequence; (2) for $n\ge 15,$ if $n$ is congruent to $3$ modulo $4$ then $T_n$ is status injective
 and among any four consecutive even orders, there is at least one order $n$ such that $T_n$ is status injective.
\end{abstract}


\section{Introduction}

We consider finite simple graphs. The {\it order} of a graph is the number of its vertices.
A connected graph is said to be {\it unicyclic} if it has exactly one cycle.
We denote  by $V(G)$ and $E(G)$ the vertex set and edge set of a graph $G$ respectively. The distance between two vertices $x$ and $y$ in a graph
is denoted by $d(x,y).$ The {\it status} of a vertex $x$ in a graph $G$, denoted by $s(x),$ is the sum of the distances between $x$ and
all other vertices; i.e.,
$$
s(x)=\sum_{y\in V(G)}d(x,y).
$$
The {\it status sequence} of $G$ is the list of the statuses of all vertices of $G$ arranged in nondecreasing order. $G$ is called {\it status injective}
 if all the statuses of its vertices are distinct [2, p.185]. Harary [4] investigated the digraph version of the concept of status in a sociometric framework,
 while Entringer, Jackson and Snyder [3] studied basic properties of this concept for graphs.

 A natural question is: Which graphs are determined by their status sequences? Slater [7] constructed infinitely many pairs of non-isomorphic trees
 with the same status sequence. Shang [5] gave a method for constructing general non-isomorphic graphs with the same status sequence.
 Let $G$ be a member of a family of graphs $\mathscr{F}$ and let the status sequence of $G$ be $s.$  $G$ is said to be {\it status unique} in
 $\mathscr{F}$ if $G$  is the unique graph in $\mathscr{F}$ whose status sequence is $s.$ Here we view two isomorphic graphs as the same graph. 
 It is known that [6] spiders are status unique in trees and that [1] status injective trees are status unique in trees.

 Shang and Lin [6, p.791] posed the following two conjectures in 2011.

{\bf Conjecture 1.} A tree and a nontree graph cannot have the same status sequence.

{\bf Conjecture 2.} Any status injective tree is status unique in all connected graphs.

In this paper we settle these two conjectures negatively. For every integer  $n\ge 10,$ we construct a tree $T_n$ and a unicyclic graph $U_n,$
both of order $n,$ with the same status sequence. There are infinitely many odd orders $n$ and infinitely many even orders $n$ such that $T_n$
is status injective.

\section{Main Results}

We will need the following lemmas. For a set $S,$ the notation $|S|$ denotes the cardinality of $S.$

{\bf Lemma 1.} [3, p.284] {\it Let $xy$ be an edge of a tree $T$ and let $T_1$ and $T_2$ be the two components of $T-xy$ with $x\in V(T_1)$
and $y\in V(T_2).$ Then $s(y)=s(x)+|V(T_1)|-|V(T_2)|.$}

{\bf Lemma 2.} {\it Let $x_0x_1x_2...x_k$ be a path in a tree and denote $d=s(x_1)-s(x_0).$ Then $s(x_{j+1})-s(x_j)\ge d+2j$ for each
$j=1,2,...,k-1.$ Consequently if $s(x_0)\le s(x_1)$ then $s(x_{j+1})-s(x_j)\ge 2j$ for each $j=1,2,...,k-1$ and in particular, $s(x_1)<s(x_2)<s(x_3)<\cdots<s(x_k).$}

{\bf Proof.} It suffices to prove the first assertion. We first show the following

Claim. If $xyz$ is a path in a tree and denote $c=s(y)-s(x),$ then $s(z)-s(y)\ge c+2.$

Let $T$ be the tree of order $n.$ Let $A$ and $B$ be the two components of $T-xy$ with $x\in V(A)$ and $y\in V(B)$,
and let $G$ and $H$ be the two components of $T-yz$ with $y\in V(G)$ and $z\in V(H).$ By Lemma 1,
$s(y)-s(x)=|V(A)|-|V(B)|=c.$ We also have $|V(A)|+|V(B)|=n$ since every edge in a tree is a cut-edge. Hence
$2|V(A)|=c+n.$ Since $V(A)\subset V(G)$ and $y\in V(G)$ but $y\notin V(A),$ we have $|V(G)|\ge |V(A)|+1.$ By Lemma 1 and
the relation $|V(G)|+|V(H)|=n$ we deduce
$$
s(z)-s(y)=|V(G)|-|V(H)|=2|V(G)|-n\ge 2|V(A)|+2-n=c+2.
$$
This proves the claim.

Applying the claim successively to the path $x_{i-1}x_ix_{i+1}$ for $i=1,2,\ldots,k-1$ we obtain the first assertion in Lemma 2.$\Box$

Lemma 2 is a generalization and strengthening of a result in [3, p.291], which states that if $x_0x_1...x_k$ is a path in a tree  and
$x_0$ has the minimum status of all vertices, then $s(x_1)<s(x_2)<\cdots < s(x_k).$

{\bf Lemma 3.} {\it The quadratic polynomial equation
$$
p^2+5p+4 = q^2+q-6
$$
in $p$ and $q$  has no nonnegative integer solution.}

{\bf Proof.} Suppose that $p$ and $q$ are nonnegative integers.
If $q\le p+2,$ then $q^2+q-6\le (p+2)^2+(p+2)-6=p^2+5p<p^2+5p+4.$ If $q\ge p+3,$ then
$q^2+q-6\ge (p+3)^2+(p+3)-6=p^2+7p+6>p^2+5p+4.$
Hence the equation cannot have any nonnegative integer solution.$\Box$

{\bf Remark.} It is not hard to prove that the only integer solutions of the equation in Lemma 3 are
$(p,q)=(-4,-3),\, (-4,2),\, (-1,-3),\, (-1,2).$

Denote by ${\mathbb N}$ the set of positive integers.

{\bf Lemma 4.} {\it Let the two functions $f(p)=p^2+5p+4$ and $h(q)=q^2+q-6$ be defined on the set ${\mathbb N}.$
If $p\ge 7$ and $|f(p)-h(q)|\le 15,$ then $q=p+2$ and $f(p)-h(q)=4.$}

{\bf Proof.} If $q\ge p+3,$ then
$$
h(q)\ge h(p+3)=f(p)+2p+2\ge f(p)+16.
$$
If $p-2\le q\le p+1,$ then
$$
f(p)\ge f(q-1)=h(q)+2q+6\ge h(q)+2p+2\ge h(q)+16.
$$
If $q\le p-3,$ then
$$
f(p)\ge f(q+3)=h(q)+10q+34\ge h(q)+44.
$$
Hence we must have $q=p+2$ and in this case, $f(p)-h(q)=4.$$\Box$

Now we are ready to state and prove the main result.

{\bf Theorem 5.} {\it For every integer $n\ge 10,$ there exists a tree $T_n$ and a unicyclic graph $U_n,$ both of order $n,$
with the following two properties:
\newline\indent (1) $T_n$ and $U_n$ have the same status sequence;
\newline\indent (2) for $n\ge 15,$ if $n\equiv 3\,\,({\rm mod}\,\, 4)$ then $T_n$ is status injective  and among any four consecutive even orders,
there is at least one order $n$ such that $T_n$ is status injective.}

{\bf Proof.} For the orders $n\ge 19$ we have a uniform construction of $T_n$ and $U_n$, and we treat this case first. For the orders $10\le n\le 18,$
the graphs will be constructed individually and they appear at the end of this proof.

Now suppose $n\ge 19.$ We distinguish the odd orders and the even orders. Let $n=2k+5$ with $k\ge 7.$ We define $T_n$ and $U_n$ as follows.
$V(T_n)=\{x_i|\,i=1,2,...,2k+5\}$ and $E(T_n)=$
$$
\{x_ix_{i+1}|\,i=1,2,...,2k-1\}\cup \{x_3x_{2k+1},x_{k+1}x_{2k+5},x_{k+3}x_{2k+4},x_{2k-3}x_{2k+3},x_{2k-2}x_{2k+2}\}.
$$
$V(U_n)=\{y_i|\,i=1,2,...,2k+5\}$ and $E(U_n)=\{y_iy_{i+1}|\,i=1,2,...,2k-1\}\cup\{y_5y_{2k+3},y_{k-1}y_{2k+4},$ $y_{k+1}y_{2k+5},y_{2k-2}y_{2k+2},y_{2k-1}y_{2k+1},y_{2k+1}y_{2k+2}\}.$
Note that $T_n$ is a caterpillar of maximum degree $3$ and $U_n$ is a unicyclic graph. $T_n$ and $U_n$ are illustrated in Figure 1.
\vskip 3mm
\par
 \centerline{\includegraphics[width=5.6in]{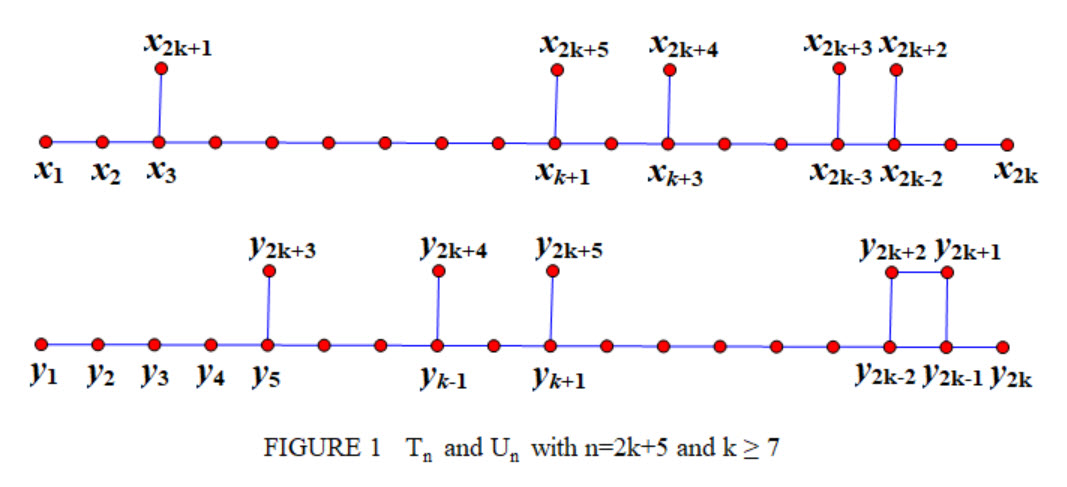}}
\par
It can be checked directly that $s(x_i)=s(y_i)$ for $i=1,2,3,k+1,2k-1,2k,...,2k+5$ and $s(x_i)=s(y_{2k+2-i})$ for $4\le i\le 2k-2.$
Hence, $T_n$ and $U_n$ have the same status sequence. For the even orders $n=2k+6$ with $k\ge 7,$ $T_n$ is obtained from $T_{n-1}$ defined above by
adding the edge $x_{2k+5}x_{2k+6},$ and $U_n$ is obtained from $U_{n-1}$ defined above by adding the edge $y_{2k+5}y_{2k+6}.$ $T_n$ and $U_n$ are
illustrated in Figure 2.
\vskip 3mm
\par
 \centerline{\includegraphics[width=5.6in]{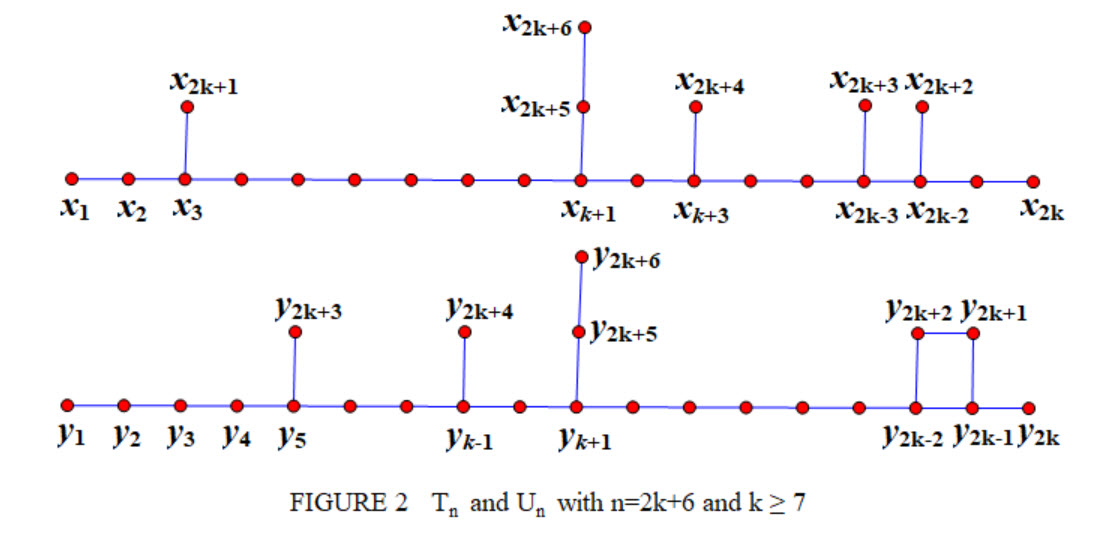}}
\par
We check easily that $s(x_i)=s(y_i)$ for $i=1,2,3,k+1,2k-1,2k,...,2k+6$ and $s(x_i)=s(y_{2k+2-i})$ for $4\le i\le 2k-2.$
Thus $T_n$ and $U_n$ also have the same status sequence.

Next we prove that the trees $T_n$ satisfy condition (2) in Theorem 5. In fact, we will determine precisely for which orders $n,$ $T_n$ is status injective.

First consider the case when $n$ is odd and let $n=2k+5$ with $k\ge 7.$ Denote $a=s(x_{k+1})=k^2+3k-2.$ We have
$$
s(x_{k-p})=\begin{cases} a+(p+2)^2-1 &{\rm if}\,\,\,\, 0\le p\le k-3,\\
                          a+k^2+1    &{\rm if}\,\,\,\, p=k-2,\\
                          a+(k+1)^2+3 &{\rm if}\,\,\,\, p=k-1;
          \end{cases}
$$
$$
s(x_{k+q})=\begin{cases} a+1 &{\rm if}\,\,\,\, q=2,\\
                          a+q^2-5    &{\rm if}\,\,\,\, 3\le q\le k-3,\\
                          a+(q+2)^2-4k+1 &{\rm if}\,\,\,\, k-2\le q\le k;
          \end{cases}
$$
$$
s(x_{2k+r})=\begin{cases} a+(k+1-r)^2+3 &{\rm if}\,\,\,\, 1\le r\le 3,\\
                          a+2k+7    &{\rm if}\,\,\,\, r=4,\\
                          a+2k+3 &{\rm if}\,\,\,\, r=5.
          \end{cases}
$$
In calculating the values $s(x_i)$ for $1\le i\le 2k$ we have used the fact that if $P=z_1z_2...z_m$ is  a path, then
$$
s(z_i)=i(i-m-1)+m(m+1)/2
$$
in $P,$ while in calculating the values $s(x_j)$ for $j=2k+1,...,2k+5$ we have used Lemma 1. From the above expressions it
follows that $x_{k+1}$ is the unique vertex with the minimum status, $x_1,x_2,x_3,x_{2k-1},x_{2k},x_{2k+1},x_{2k+2},x_{2k+3}$
are the vertices with the eight largest statuses, since
$$
s(x_1)>s(x_{2k})>s(x_{2k+1})>s(x_2)>s(x_{2k+2})>s(x_{2k-1})>s(x_3)>s(x_{2k+3})>s(x_i) \eqno (1)
$$
for any $i\neq 1,2,3,2k-1,2k,2k+1,2k+2,2k+3$ and
$$
s(x_{2k+1})>s(x_{2k+2})>s(x_{2k+3})>s(x_{2k+4})>s(x_{2k+5}). \eqno (2)
$$
Partition the vertex set of $T_n$ into three sets:
$$
L=\{x_i|\, 1\le i\le k\},\,\,\, R=\{x_i|\, k+1\le i\le 2k\}\,\,\, {\rm and}\,\,\,W=\{x_i|\, 2k+1\le i\le 2k+5\}.
$$
The inequalities in (2) show that any two distinct vertices in $W$ have different statuses. Applying Lemma 2 to the two paths
$x_{k+1}x_kx_{k-1}...x_2x_1$ and $x_{k+1}x_{k+2}...x_{2k-1}x_{2k}$ we see that any two distinct vertices in $L$ or in $R$ have different statuses.
Next we show that for any $x\in L$ and $y\in R,$ $s(x)\neq s(y).$ By the inequalities in (1) it suffices to prove that
$s(x_i)\neq s(x_j)$ for $4\le i\le k$ and $k+2\le j\le 2k-2,$ which is equivalent to $s(x_{k-p})\neq s(x_{k+q})$ for $0\le p\le k-4$ and
$2\le q\le k-2.$ We have the expressions $s(x_{k-p})=a+(p+2)^2-1$ for $0\le p\le k-4,$ $s(x_{k+2})=a+1,$ $s(x_{2k-2})=a+k^2-4k+1$
and $s(x_{k+q})=a+q^2-5$ for $3\le q\le k-3.$ First, $s(x_{k-p})\ge s(x_k)=a+3>a+1=s(x_{k+2}).$
The equality $s(x_{k-p})=s(x_{k+q})$ for $3\le q\le k-3$ is equivalent to $4=(q+p+2)(q-p-2),$ which is impossible, since $q+p+2\ge 5$ and $q-p-2$ is an integer. Also, $s(x_{k-p})=s(x_{2k-2})$ is equivalent to $2=(k+p)(k-p-4),$ which is impossible, since $k+p\ge 7$ and $k-p-4$ is an integer.
Hence $s(x)\neq s(y)$ for $x\in L$ and $y\in R.$

By the above analysis, it is clear that the only possibilities for two distinct vertices to have the same status are $s(x_{2k+5})=s(x_i)$ and
$s(x_{2k+4})=s(x_i)$ for $4\le i\le k$ or $k+2\le i\le 2k-2.$ By the expressions for their status values, it is easy to verify that
$s(x_{2k+5})=s(x_i)$ for some $i$ with $4\le i\le k$ if and only if $k=2c^2-2$ for some integer $c;$
$s(x_{k+2})<s(x_{2k+5})<s(x_{2k-2})$ and $s(x_{2k+5})=s(x_i)$ for some $i$
with $k+3\le i\le 2k-3$ if and only if $k=2c^2-4$ for some integer $c;$ $s(x_{2k+4})=s(x_i)$ for some $i$ with $4\le i\le k$ if and
only if $k=2c^2-4$ for some integer $c;$ $s(x_{k+2})<s(x_{2k+4})<s(x_{2k-2})$
and $s(x_{2k+4})=s(x_i)$ for some $i$ with $k+3\le i\le 2k-3 $ if and only if $k=2c^2-6$ for some integer $c.$

Thus, $T_n$ with $n=2k+5$ is not status injective if and only if $k=2c^2-2,$  $2c^2-4$ or $2c^2-6$ for some integer $c.$ Since all these values
of $k$ are even, it follows that for every odd $k,$ $T_n$ is status injective; i.e., if $n\equiv 3\,\,({\rm mod}\,\, 4)$ then $T_n$
is status injective.

Next we treat the case when the order $n$ is even. Let $n=2k+6$ with $k\ge 7.$ With $d=s(x_{k+1})=k^2+3k$ we have
$$
s(x_{k-p})=\begin{cases} d+p^2+5p+4 &{\rm if}\,\,\,\, 0\le p\le k-3,\\
                          d+k^2+k    &{\rm if}\,\,\,\, p=k-2,\\
                          d+k^2+3k+4 &{\rm if}\,\,\,\, p=k-1;
          \end{cases}
$$
$$
s(x_{k+q})=\begin{cases} d+2 &{\rm if}\,\,\,\, q=2,\\
                          d+q^2+q-6    &{\rm if}\,\,\,\, 3\le q\le k-3,\\
                          d+q^2+5q-4k+4 &{\rm if}\,\,\,\, k-2\le q\le k;
          \end{cases}
$$
$$
s(x_{2k+r})=\begin{cases} d+k^2+k+2 &{\rm if}\,\,\,\, r=1,\\
                          d+k^2-k+2    &{\rm if}\,\,\,\, r=2,\\
                          d+k^2-3k+4    &{\rm if}\,\,\,\, r=3,\\
                          d+2k+10    &{\rm if}\,\,\,\, r=4,\\
                          d+2k+2 &{\rm if}\,\,\,\, r=5,\\
                          d+4k+6    &{\rm if}\,\,\,\, r=6.\\
          \end{cases}
$$
From the above expressions we deduce that $x_{k+1}$ is the unique vertex with the minimum status $d.$ The case $k=7$ corresponds to
$n=20$ and we check directly that $T_{20}$ is status injective. Next suppose $k\ge 8.$ Then
$x_1,x_2,x_3,x_{2k-1},x_{2k},x_{2k+1},x_{2k+2},x_{2k+3}$ are the vertices with the eight largest statuses, since
$$
s(x_1)>s(x_{2k})>s(x_{2k+1})>s(x_2)>s(x_{2k+2})>s(x_{2k-1})>s(x_3)>s(x_{2k+3})>s(x_i) \eqno (3)
$$
for any $i\neq 1,2,3,2k-1,2k,2k+1,2k+2,2k+3.$ Also
$$
s(x_{2k+1})>s(x_{2k+2})>s(x_{2k+3})>s(x_{2k+6})>s(x_{2k+4})>s(x_{2k+5}). \eqno (4)
$$
In considering two vertices with equal status, we can exclude the eight vertices with the eight largest statuses by (3) and the unique vertex
$x_{k+1}$ with the minimum status. Denote
$$
L^{\prime}=\{x_i|\, 4\le i\le k\},\,\,\, R^{\prime}=\{x_i|\, k+2\le i\le 2k-2\}\,\,\, {\rm and}\,\,\,W^{\prime}=\{x_i|\, 2k+1\le i\le 2k+6\}.
$$
Let $x$ and $y$ be two distinct vertices with $s(x)=s(y).$ By the inequalities in (4), it is impossible that $x,y\in W^{\prime}.$
By Lemma 2 we cannot have $x,y\in L^{\prime}$ or $x,y\in  R^{\prime}.$ Suppose $x\in L^{\prime}$ and $y\in R^{\prime}.$ We have $s(x)>s(x_{k+2}),$ $s(x_4)>s(x_{2k-2})$ and $s(x_i)<s(x_{2k-2})$ for $5\le i\le k.$ Thus, $y\neq x_{k+2},\,x_{2k-2}.$ We have $x=x_i$ for some $i$ with $4\le i\le k$ and $y=x_j$
for some $j$ with $k+3\le j\le 2k-3.$ Hence $s(x)=d+p^2+5p+4$ with $0\le p\le k-4$ and $s(y)=d+q^2+q-6$ with $3\le q\le k-3.$ Then $s(x)=s(y)$ yields $p^2+5p+4=q^2+q-6,$ which is impossible by Lemma 3.

Now, by (3) and the above analysis it is clear that $s(x)=s(y)$ can occur only if $x\in\{x_{2k+4},x_{2k+5},x_{2k+6}\}$ and
$y\in L^{\prime}\cup R^{\prime}$ or the roles of $x$ and $y$ are interchanged. The case $k=8$ corresponds to $n=22,$ and we check directly that
$T_{22}$ is not status injective. Next we suppose $k\ge 9.$ Then $s(x_{2k-2})>s(x_{2k+6})>s(x_{2k+4})>s(x_{2k+5}),$ and hence $x_{2k-2}$ can be excluded
from $R^{\prime}.$ Similarly, since $s(x_{k+2})<s(x_k)<s(x_{2k+5})<s(x_{2k+4})<s(x_{2k+6}),$  $x_k$ can be excluded from $L^{\prime}$
and $x_{k+2}$ can be excluded from $R^{\prime}.$
Note that the statuses of the vertices in $L^{\prime}\setminus \{x_k\}$ have the uniform expression $d+p^2+5p+4$ with $1\le p\le k-4$
and the statuses of the vertices in $R^{\prime}\setminus \{x_{k+2},\,x_{2k-2}\}$ have the uniform expression $d+q^2+q-6$
with $3\le q\le k-3.$

Denote the empty set by $\phi,$ and denote
$\Omega_k=\{2k+2,\, 2k+10,\, 4k+6\},$ $\Gamma_k=A_k\cup B_k$ where $A_k=\{p^2+5p+4|\,1\le p\le k-4,\, p\in {\mathbb N}\}$ and
$B_k=\{q^2+q-6|\,3\le q\le k-3,\,q\in {\mathbb N}\}.$ It follows that when $k\ge 9,$ $T_n$ has two distinct vertices with the same status
if and only if $\Omega_k \cap \Gamma_k \neq\phi.$ Denote $\Gamma=A\cup B$ where $A=\{p^2+5p+4|\,p\in {\mathbb N}\}$ and $B=\{q^2+q-6|\,q\in {\mathbb N}\}.$
Since $\Omega_k \cap \Gamma_k=\Omega_k \cap \Gamma,$ we obtain the following criterion for $k\ge9:$

$T_n$ is status injective if and only if $\Omega_k \cap \Gamma =\phi.$

The graphs $T_n$ with $15\le n\le 18$ constructed below are all status injective. Using the above criterion we can check that $T_n$ is status injective
for
$$
k=10,14,18,21,23,25,27,29,33,35,38,40,42.
$$
Thus the assertion in Theorem 5 on $T_n$ for even $n$ with $k\le 42$ is true.

Next we suppose $k\ge 43.$ We will prove that among the four numbers $k,k+1,k+2,k+3$ there is at least one for which $T_n$ is status injective.
To do so, consider
\begin{eqnarray*}
\Omega_k&=&\!\{2k+2,\,2k+10,\,4k+6\}\\
\Omega_{k+1} &=& \{2k+4,\,2k+12,\,4k+10\}\\
\Omega_{k+2} &=& \{2k+6,\,2k+14,\,4k+14\}\\
\Omega_{k+3} &=& \{2k+8,\,2k+16,\,4k+18\}.
\end{eqnarray*}
 The numbers in these four sets can be
partitioned into two classes:
$$
X=\{2k+i|\,i=2,4,6,8,10,12,14,16\}\,\,\,\,{\rm and}\,\,\,\, Y=\{4k+j|\,j=6,10,14,18\}.
$$
We claim that
$$
|X\cap A|\le 1,\,\,\,\,|X\cap B|\le 1,\,\,\,\,|Y\cap A|\le 1,\,\,\,\,|Y\cap B|\le 1. \eqno (5)
$$
Define two polynomials $f(p)=p^2+5p+4$ and $h(q)=q^2+q-6.$ Then $A=\{f(p)|\,p\in {\mathbb N}\}$ and $B=\{h(q)|\,q\in {\mathbb N}\}.$
In the sequel the symbol $\Rightarrow$ means ``implies".
We first prove $|X\cap A|\le 1.$ To the contrary, suppose there exist $i,j,p_1,p_2$ with $2\le i<j\le 16$ and $p_1<p_2$ such that
$f(p_1)=2k+i$ and $f(p_2)=2k+j.$ $k\ge 43$ and $i\ge 2$ $\Rightarrow$ $f(p_1)=2k+i\ge 88$ $\Rightarrow$ $p_1\ge 7.$
We have $f(p_2)-f(p_1)=j-i\le 14.$ But on the other hand, $f(p_2)-f(p_1)\ge f(p_1+1)-f(p_1)= 2p_1+6\ge 20,$ a contradiction.
The inequality $|X\cap B|\le 1$ is similarly proved by using the fact that $h(q)\in X$ $\Rightarrow$ $h(q)\ge 88$ $\Rightarrow$ $q\ge 10.$
The inequalities $|Y\cap A|\le 1$ and $|Y\cap B|\le 1$ can also be similarly proved by using the facts that
$f(p)\in Y$ $\Rightarrow$ $f(p)\ge 178$ $\Rightarrow$ $p\ge 11$ and $h(q)\in Y$ $\Rightarrow$ $h(q)\ge 178$ $\Rightarrow$ $q\ge 14.$

Note that the assumption $k\ge 43$ implies that ${\rm min}\,X\ge 88$ and ${\rm min}\,Y\ge 178.$ Hence if $f(p)\in X\cup Y$ we have
$p\ge 7$ and Lemma 4 can be applied.

Suppose $\Omega_i\cap\Gamma\neq\phi$ for $i=k,k+1,k+2.$ We will show that $\Omega_{k+3}\cap\Gamma =\phi.$ Since $\Omega_k\cap \Gamma\neq\phi,$
at least one of the two cases $\{2k+2,\,2k+10\}\cap \Gamma\neq\phi$ and $4k+6\in \Gamma$ must occur. Recall that $\Gamma=A\cup B.$

Case 1. $\{2k+2,\,2k+10\}\cap \Gamma\neq\phi.$ We first consider the case when $\{2k+2,\,2k+10\}\cap A\neq\phi.$
Denote $\Psi=\{2k+4,\, 2k+12,\, 2k+8,\, 2k+16\}.$ By (5), $\Psi\cap A=\phi.$ By Lemma 4, $\Psi\cap B=\phi.$
It follows that $\Psi\cap\Gamma=\phi.$ Since $\Omega_{k+1}\cap \Gamma\neq\phi$ and $\{2k+4,\, 2k+12\}\cap\Gamma=\phi,$
we deduce that $4k+10\in\Gamma.$ By (5), $4k+10$ and $4k+18$ can not be both in $A$ or both in $B.$ Since $4\neq 8=(4k+18)-(4k+10)\le 15,$
by Lemma 4 it is also impossible that one of $4k+10$ and $4k+18$ is in $A$ and the other in $B.$ But $4k+10\in\Gamma=A\cup B.$
Hence $4k+18\notin\Gamma$ and we obtain $\Omega_{k+3}\cap\Gamma=\phi.$ The case when $\{2k+2,\,2k+10\}\cap B\neq\phi$ is similar. Again
we use (5), Lemma 4 and $\Omega_{k+1}\cap\Gamma\neq\phi$ to deduce $\Omega_{k+3}\cap\Gamma=\phi.$

Case 2. $4k+6\in \Gamma.$ Using (5) and Lemma 4 we deduce that $\{4k+14,\, 4k+18\}\cap\Gamma=\phi.$ Then the condition $\Omega_{k+2}\cap\Gamma\neq\phi$
implies $\{2k+6,\,2k+14\}\cap\Gamma\neq\phi.$ Applying (5) and Lemma 4 once more
we have $\{2k+8,\,2k+16\}\cap\Gamma=\phi.$ Hence $\Omega_{k+3}\cap\Gamma=\phi.$

This completes the proof of the case $n\ge 19$ of Theorem 5. The graph pairs $T_n$ and $U_n$ with $10\le n\le 18$ are depicted in
Figures 3-11 below. They satisfy the condition $s(T_n)=s(U_n)$ and for $15\le n\le 18,$ $T_n$ is status injective. In these graphs,
the number beside a vertex is the status of that vertex.
\vskip 3mm
\par
 \centerline{\includegraphics[width=5in]{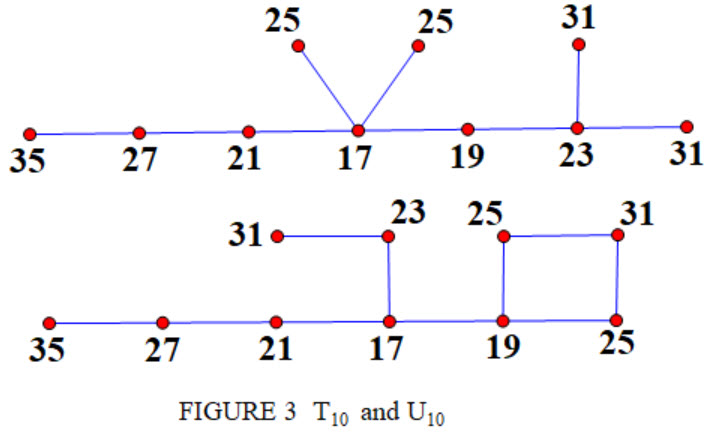}}
\par
\vskip 3mm
\par
 \centerline{\includegraphics[width=5in]{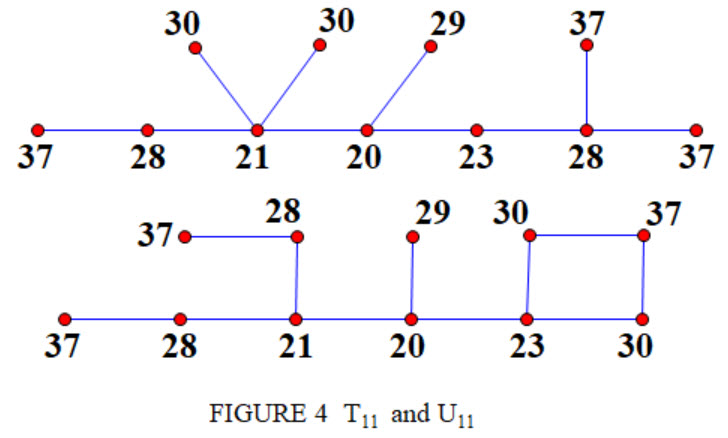}}
\par
\vskip 3mm
\par
 \centerline{\includegraphics[width=5in]{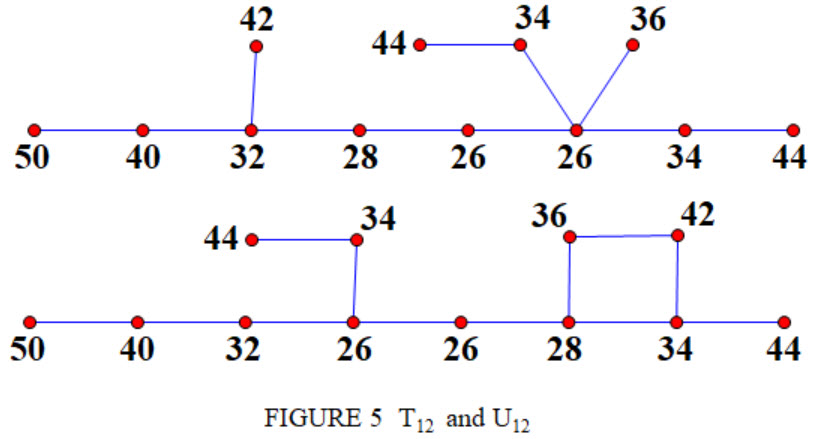}}
\par
\vskip 3mm
\par
 \centerline{\includegraphics[width=5in]{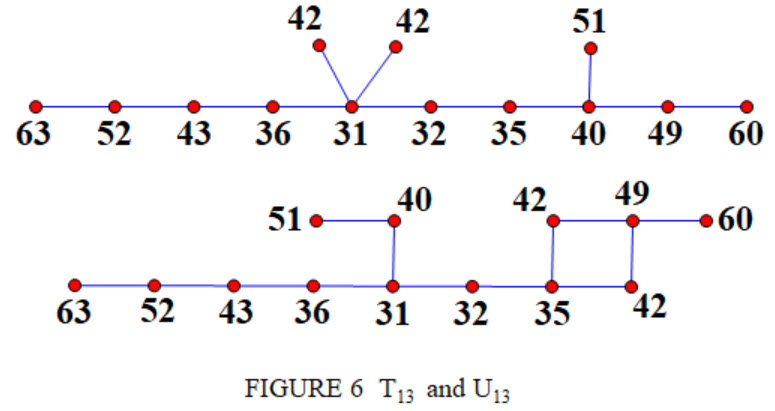}}
\par
\vskip 4mm
\par
 \centerline{\includegraphics[width=5in]{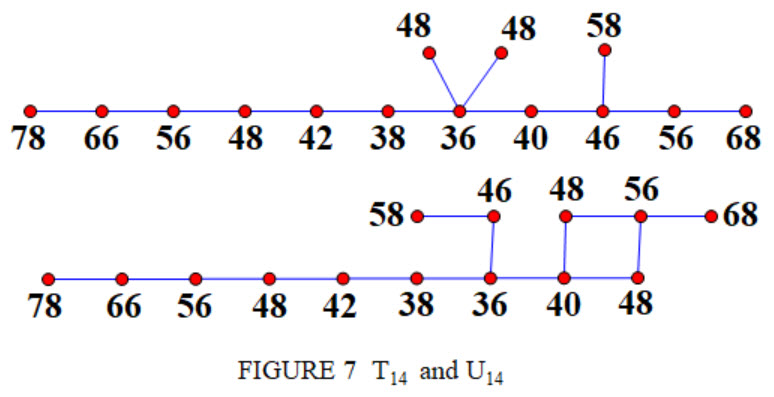}}
\par
\vskip 4mm
\par
 \centerline{\includegraphics[width=5in]{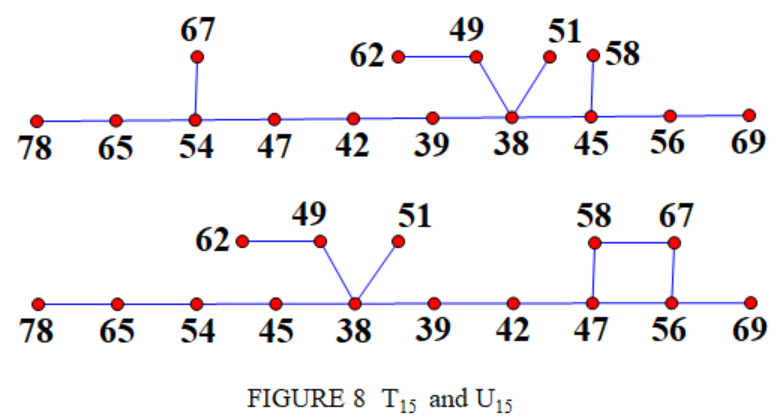}}
\par
\vskip 4mm
\par
 \centerline{\includegraphics[width=5in]{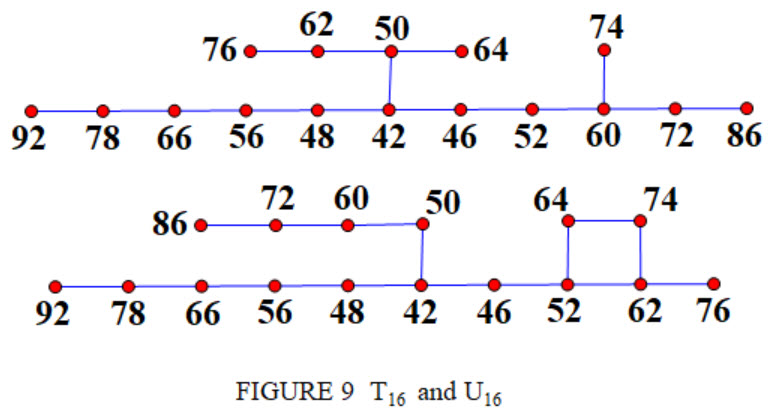}}
\par
\vskip 4mm
\par
 \centerline{\includegraphics[width=5in]{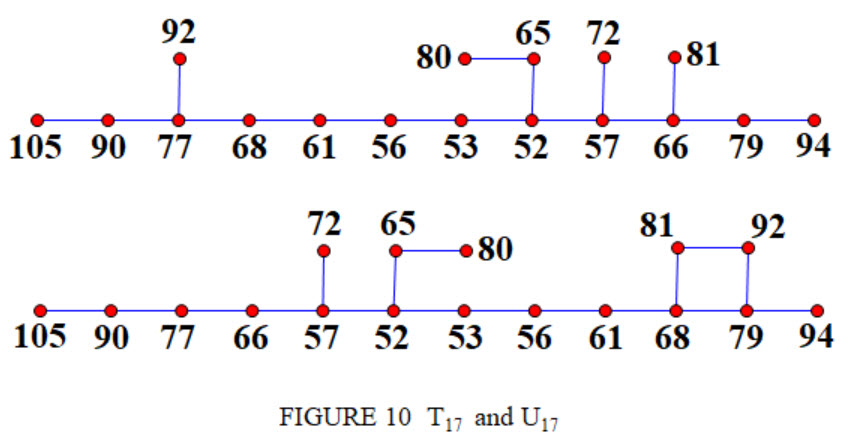}}
\par
\vskip 4mm
\par
 \centerline{\includegraphics[width=5in]{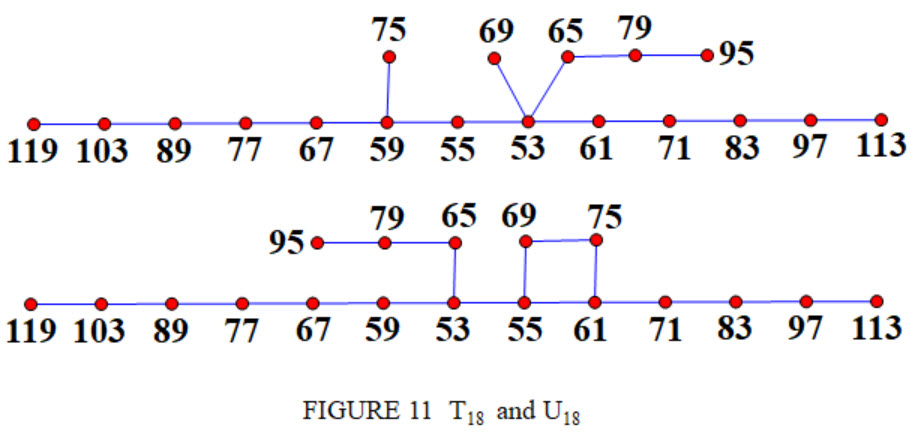}}
\par
This completes the proof of Theorem 5.$\Box$

{\bf Remark.} A computer search shows that $10$ is the smallest order for the existence of a tree and a nontree graph with
the same status sequence.

\vskip 5mm
{\bf Acknowledgement.} This research  was supported by the NSFC grants 11671148 and 11771148 and Science and Technology Commission of Shanghai Municipality (STCSM) grant 18dz2271000.

\end{document}